\documentclass[12pt]{article}
\usepackage{a4,amsmath,amsthm,amssymb,amsfonts,enumerate,hyperref}

\theoremstyle{plain}

\newtheorem{thm}{Theorem}      
\newtheorem{cor}{Corollary}
\newtheorem{lem}{Lemma}  
\newtheorem{conj}{Conjecture}

\theoremstyle{definition}

\newcommand{\abs}[1]{{\left| {#1} \right|}}
\newcommand{\p}[1]{{\left( {#1} \right)}}

\author{Johan Andersson\thanks{Department of Mathematics, Stockholm University, SE-10691, Sweden. {\it johana@math.su.se}}}

\begin{document}

\title{On the solutions to a power sum problem}


\maketitle

\begin{abstract}
In a recent paper \cite{Andersson2} we proved that
\begin{gather*}
 \inf_{ |z_k| = 1}
     \max_{\nu=1,\dots,n^2-n}   \left| \sum_{k=1}^n z_k^\nu  \right| = \sqrt{n-1}
\end{gather*}
  when $n-1$ is  a prime power. We also proved that a construction of Fabrykowski of an explicit $n$-tuple $(z_1,\ldots,z_n)$ provided a global minimum to the problem.  The construction depends on the construction of perfect difference sets of order $n-1$.  As an open problem we asked whether we could obtain all global minima by this construction. In this paper we will show that this is in fact true and we furthermore prove that
\begin{gather*}
 \inf_{ |z_k| = 1}
     \max_{\nu=1,\dots,n^2-n}   \left| \sum_{k=1}^n z_k^\nu  \right|> \sqrt{n-1}
\end{gather*}
in case there exist no perfect difference set of order $n-1$. A well 
 known conjecture states that each perfect difference set must have prime power order.
\end{abstract}

\section{Introduction}
Throughout this paper we will let 
\begin{gather} S(\nu)= \sum_{k=1}^n z_k^\nu \end{gather} denote the pure power sums. 
In the papers \cite{Andersson1},\cite{Andersson2} and \cite{Andersson3} we considered Tur\'an's problem 10 which is to estimate
\begin{gather*}
 (*)=\inf_{ |z_k| \geq 1}
     \max_{\nu=1,\dots,n^2}   \abs{S(\nu)}.
\end{gather*}
We showed already in \cite{Andersson1} that  $\sqrt n \leq (*) \leq \sqrt{n+1}$ if $n+1$ is a prime. In general we have that $\sqrt n \leq (*)$ 
(see \cite{Andersson1}) and $(*)=\sqrt n +O(n^{0.2625+\epsilon})$ (see \cite{Andersson3}).  
It would be interesting to find an explicit solution (in terms of minimal systems $(z_1,\ldots,z_n)$) to this inf max problem, but due to lack of success in this we considered (see \cite{Andersson2}) the modified problems
\begin{gather*}
 (**)=\inf_{ |z_k|  \geq 1}
     \max_{\nu=1,\dots,n^2-2}   \left| S(\nu) \right|, \\ \intertext{and}
   (***)= \inf_{ |z_k|  = 1}
     \max_{\nu=1,\dots,n^2-n}   \left| S(\nu) \right|.
\end{gather*} 
We proved that $(**)=\sqrt{n}$ if $n$ is a prime power and $(***)=\sqrt{n-1}$ if $n-1$ is a prime power. We
 gave explicit solutions and asked whether these were the only solutions possible. In this paper we will answer this question in the case $(***)$ (\cite{Andersson2}, Problem 2). In our paper we used a theorem independently proved  by Newman, Cassels and Szalay (Szalay \cite{Szalay} and Theorem 7.3 of Tur\'an \cite{Turan})
to get a lower bound, and our solution will come from examining  the steps of the proof of that theorem (compare the proof of Lemma \ref{ett} below and page 81 in Tur\'an \cite{Turan}) and modifying the proof to better fit our purpose.

\section{A lower bound}

\begin{lem} \label{ett}
 Suppose that $(z_1,\ldots,z_n)$ is an $n$-tuple of unimodular complex numbers. One then has that
 \begin{enumerate}[(i)]
  \item  $\max_{\nu=1,\dots,n^2-n}   \left| S(\nu) \right| \geq \sqrt {n-1}$.
  \item If the inequality in $(i)$ holds with equality, then $\abs{S(\nu)}=\sqrt{n-1}$ for all $\nu=1,\ldots,n^2-n$.
  \end{enumerate}
\end{lem}
\begin{proof}
Let $\theta_k$ be defined so that $z_k=e(\theta_k)$. Then  
\begin{gather*}
  S(\nu)=\sum_{k=1}^n e(\nu \theta_k). \qquad \qquad (e(x)=e^{2 \pi i x})
\end{gather*}
Let
 \begin{gather*}
   F_m(t)=\sum_{\nu=1-m}^{m-1}\p{1-\frac{\abs{\nu}}m} e(\nu t)
\end{gather*}
be the $m$'th Fej\'er kernel. The Fej\'er kernel can be written as
\begin{gather} \label{nonneg}
   F_m(t) =\frac 1 {m} \left( \frac{\sin \pi mx}{\sin \pi x} \right)^2,
\end{gather} 
and  is thus non negative. We have that
\begin{gather*} 
  \sum_{\nu=1-m}^{m-1}  \p{1-\frac{\abs{\nu}}m}  \abs{S(\nu)}^2 = \sum_{k,l=1}^n F_m(\theta_{k}-\theta_{l}),  \\ 
\intertext{which by the contribution of the diagonal $k=l$, and the non negativity of  the Fej\'er kernel implies that}
 \sum_{\nu=1-m}^{m-1}  \p{1-\frac{\abs{\nu}}m} \abs{S(\nu)}^2 \geq n F_m(0).
\end{gather*}
By  letting $m=n^2-n+1$ and  $\abs{S(\nu)}^2 = n-1+\varepsilon_\nu$ we get
\begin{gather*}
  \sum_{\nu=n-n^2}^{n^2-n} \p{1-\frac{\abs{\nu}}{n^2-n+1}}  \varepsilon_\nu \geq F_{n^2-n+1}(0).
\end{gather*}
By the fact that $\abs{S(\nu)}=\abs{S(\nu)}$ we have that  $\varepsilon_{-\nu}=\varepsilon_\nu$. By furthermore using the facts that  
$F_{n^2-n+1}(0)=n^2-n+1$ and $\varepsilon_0=n^2-n+1$, we obtain the following inequality
\begin{gather*}
  \sum_{\nu=1}^{n^2-n}  \p{1-\frac{\nu}{n^2-n+1}}  \varepsilon_\nu \geq 0.
\end{gather*}
This implies that for at least one $1 \leq \nu \leq n^2-n$ we have that $\varepsilon_\nu \geq 0$,  which implies $(i)$. 
Furthermore if all $\varepsilon_\nu \leq 0$ then $\varepsilon_\nu=0$ for all $\nu=1,\ldots,k$ which implies $(ii)$
\end{proof}

\begin{lem} \label{ngon} 
   Let $(z_1,\ldots,z_n)$ be an $n$-tuple of unimodular complex numbers. Then the following statements are equivalent:
 \begin{enumerate}[(i)] \item $S(\nu)=0$ for $\nu=1,\ldots,n-1$.
   \item  The numbers $z_k$ are the $n$ vertices of a regular $n$-gon on the unit circle.
\end{enumerate}
\end{lem}

\begin{proof}
  The fact that $(ii) \implies (i)$ is immediate. Remains to prove that $(i) \implies (ii)$.  Consider the polynomial
  \begin{gather*}
    p(x)=\prod_{k=1}^n (x-z_k) = x^n+a_1 x^{n-1}+ \cdots+a_n.
  \end{gather*}
   By the Newton-Girard identities
\begin{gather*}
  S(\nu)+a_1 S(\nu-1)+\dots+a_{\nu-1} S(1)+\nu a_\nu=0, \qquad \qquad (\nu=1,2,\dots,n) \end{gather*}
 and the fact that $S(\nu)=0$ for $\nu=1,\ldots,n-1$ we find that $a_\nu=0$ for $\nu=1,\dots,n-1$. We see that
\begin{gather*}   
p(x)=x^n+(-z_1) \cdots (-z_n)= x^n-w
 \end{gather*}
 for some complex number $w$ with $|w|=1$.
\end{proof}

\section{Perfect Difference sets}

\begin{lem} {\rm (Singer \cite{Singer})} \label{Singer} Let $q$ be a prime power. Then there exists a perfect difference set of order $q$ or in other words
integers $a_1,\ldots,a_{q+1}$
such that the integers $a_i-a_j$ for $i \neq j$ form all non zero
residues $\mod q^2+q+1$.
\end{lem}

\begin{proof} This follows by the use of projective planes over finite fields (See Singer \cite{Singer}). 
\end{proof}
\begin{conj} \label{harg} {\rm (The prime power conjecture)} There only exist perfect different sets of order $q$ if $q$ is a prime power.
\end{conj}
The prime power conjecture has been proved for a lot of special cases. If $n \equiv 1,2 \pmod 4$ and can not be written as a sum of two squares then 
there exist no perfect difference set of order $n$ (Bruck-Ryser \cite{BruckRyser}). In case $n\geq 6$ and $n \equiv 3,6 \pmod 9$ (Wilbrink \cite{Wilbrink}). Finally it has been proved for $n<2 \cdot 10^9$ (Baumert-Gordon \cite{BaumertGordon}). We remark here that the problem that a general finite projective plane has prime power order is  more difficult (The Bruck-Ryser result is still true, but other methods fails), and is yet to be proved true for $n=12$. Conjecture \ref{harg} is equivalent to the conjecture that a finite  {\em cyclic} projective plane must have prime power order.

\section{Main result}

\begin{thm} \label{thmett}
 Let $(z_1,\ldots,z_n)$ be an $n$-tuple of unimodular complex numbers. Then the following statements are equivalent:
\begin{enumerate}[(i)]
 \item $\abs{S(\nu)}=\sqrt{n-1}$ for $\nu=1,\ldots,n^2-n$.
 \item There exist a complex number $\alpha$ with $|\alpha|=1$ and a perfect difference set $(a_1,\ldots,a_n)$ of order $n-1$ (which means that  $a_i-a_j$ for $i \neq j$ form all non zero residues $\mod n^2-n+1$) such that   $\displaystyle z_k=\alpha \, e\p{\frac{a_k}{n^2-n+1}}$. 
\end{enumerate}
 In case $(i)$ and $(ii)$ are false then one has that $\displaystyle \max_{\nu=1,\ldots,n^2-n} \abs{S(\nu)}>\sqrt{n-1}$.
\end{thm}

\begin{proof}
\noindent By choosing 
\begin{gather*}  z_k=\alpha e \p{\theta_k}, \qquad \qquad \p{e(x)=e^{2 \pi i x}} \\
\intertext{with}
  \theta_k=\frac{a_k}{n^2-n+1}, \qquad \theta_1=0 \qquad \text{and} \qquad |\alpha|=1, \intertext{we see that} \notag \begin{split} \abs{\sum_{k=1}^n z_k^\nu}^2 &=
n+\sum_{k \neq l} e \p{\nu(\theta_k-\theta_l)},  \\ &= n-1 +
\sum_{j=1}^{n^2-n+1} e \p{\nu \lambda_j},
\end{split} \\ \intertext{where $\lambda_1=0$ and $0  \leq \lambda_2 \leq \cdots \leq \lambda_{n^2-n+1} \leq 1$ consists of the $n^2-n$ differences $\theta_k-\theta_l$ where $k \neq l$. Then }
   \abs{S(\nu)}=\sqrt{n-1} \qquad \qquad (\nu=1,\ldots,n^2-n) \\ \intertext{is equivalent to} \sum_{j=1}^{n^2-n+1} e \p{\nu \lambda_j} =0. \qquad \qquad (\nu=1,\ldots,n^2-n)
\end{gather*} 
By Lemma \ref{ngon} this is true iff $\lambda_j$ consists of the vertices of a regular $(n^2-n+1)$-gon. Since $\lambda_0=1$ and $\lambda_j$ is increasing this is true iff $\displaystyle \lambda_j=\frac{j}{n^2-n+1}$. This is equivalent to the fact that  $a_i-a_j$ for $i \neq j$ form all non zero residues $\mod n^2-n+1$. That the $a_j$'s are furthermore integers follows by the fact that $\displaystyle \theta_k=\theta_k-0=\theta_k-\theta_0 = \frac{j_k}{n^2-n+1}$ for integers $j_k$. 
This means that $(a_1,\ldots,a_n)$ form a perfect difference set of order $n-1$. 
 In case $\abs{S(\nu)} \neq \sqrt{n-1}$ for some $1 \leq \nu \leq n^2-n$ then it follows from Lemma \ref{ett} that $\abs{S(\nu)}>\sqrt{n-1}$ for some   $1 \leq \nu \leq n^2-n$.
\end{proof}

As consequences of our theorem we have the following corollaries:
\begin{cor} Suppose that there exist no perfect difference set of order $n-1$. Then one has that
\begin{gather*}
 \inf_{ |z_k| = 1}
     \max_{\nu=1,\dots,n^2-n}   \left| S(\nu) \right| > \sqrt{n-1}. 
\end{gather*}
\end{cor}

\begin{cor} The prime power conjecture is equivalent to the following statement: Suppose that $n-1$ is not a prime power. Then
 \begin{gather*}
 \inf_{ |z_k| = 1}
     \max_{\nu=1,\dots,n^2-n}   \left| S(\nu) \right|> \sqrt{n-1}.
\end{gather*}
\end{cor}

\begin{cor} Suppose that $n-1$ is a prime power. Then one has that  
\begin{gather*}   \inf_{ |z_k| = 1}
     \max_{\nu=1,\dots,n^2-n} \left| S(\nu) \right|=\sqrt{n-1}, 
\end{gather*} 
and every minimal system is given by the construction $\displaystyle z_k=\alpha e\left(\frac{a_k} {n^2-n+1} \right)$ of Fabrykowski with
 $|\alpha|=1$ and where $(a_1,\ldots,a_n)$ is a perfect difference set of order $n-1$.
\end{cor}
\begin{proof} This follows from Theorem \ref{thmett} and Lemma \ref{Singer} with $q=n-1$. Further discussion of the Fabrykowski construction is given in \cite{Fabrykowski} and \cite{Andersson2}.\end{proof}

\providecommand{\bysame}{\leavevmode\hbox to3em{\hrulefill}\thinspace}
\providecommand{\href}[2]{#2}
\bibliographystyle{halpha}

\end{document}